\documentclass[a4paper,12pt,frenchb,english]{amsart}
\usepackage[latin1]{inputenc}
\usepackage[T1]{fontenc}
\usepackage{amssymb,lmodern,calrsfs,babel}
\usepackage{hyperref}
\NoAutoSpaceBeforeFDP

\newtheorem{thm}{Th{\'e}or{\`e}me}

\theoremstyle{remark}
\newtheorem{rem}[thm]{Remarque}

\newcommand{\II}{\mathbb{I}}
\newcommand{\LL}{\mathcal{L}}
\newcommand{\bg}{\bar{g}}
\newcommand{\bga}{\bar{\gamma}}
\newcommand{\dg}{\dot{g}}
\newcommand{\bn}{{\bar{\nabla}}}

\DeclareMathOperator{\Tr}{Tr}
\DeclareMathOperator{\Ric}{Ric}
\DeclareMathOperator{\vol}{vol}
\newcommand{\bRic}{\overline{\mathrm{Ric}}}

\renewcommand{\leq}{\leqslant}
\renewcommand{\geq}{\geqslant}

\begin{document}

\title[Continuation unique pour les m{\'e}triques d'Einstein]{Continuation
  unique {\`a} partir de l'infini conforme pour les m{\'e}triques d'Einstein}
\author{Olivier Biquard}
\address{Universit{\'e} Pierre et Marie Curie-Paris 6, UMR 7586, Institut
  de Math{\'e}matiques de Jussieu}
\thanks{L'auteur est soutenu par le contrat 06-BLAN60154-01 de l'ANR}
\maketitle

\begin{abstract}
  We prove the unique continuation property at the conformal infinity
  for asymptotically hyperbolic Einstein metrics.
\end{abstract}

\selectlanguage{frenchb}
Le but de cette note est de prouver un r{\'e}sultat d'unicit{\'e} pour les
m{\'e}triques d'Einstein asymptotiquement hyperboliques dont le
comportement {\`a} l'infini est fix{\'e} {\`a} un ordre suffisamment grand. Cette
question, ainsi que sa variante infinit{\'e}simale, semble pos{\'e}e depuis
quelque temps, et est importante dans le cadre du programme d'Anderson
\cite{And05} visant {\`a} donner des crit{\`e}res g{\'e}n{\'e}raux d'existence de
telles m{\'e}triques {\`a} bord conforme donn{\'e}.

Soit une vari{\'e}t{\'e} $M^{n+1}$ {\`a} bord $X^n=\partial M$, une m{\'e}trique $g$ dans
l'int{\'e}rieur de $M$ est dite asymptotiquement hyperbolique s'il existe
une {\'e}quation $\{x=0\}$ de $X\subset M$, et une m{\'e}trique $\gamma$ sur $X$, telles
que, pour $x\to0$,
\begin{equation}
  \label{eq:1}
  g \sim \frac{dx^2+\gamma}{x^2} .
\end{equation}
Seule la classe conforme $[\gamma]$ de $\gamma$ est alors bien d{\'e}finie
par $g$ et est appel{\'e}e l'infini conforme de $g$. Un exemple de telle
m{\'e}trique est {\'e}videmment la m{\'e}trique hyperbolique
$\frac{dx^2+dx_1^2+\cdots+dx_n^2}{x^2}$. Plus g{\'e}n{\'e}ralement, {\`a} partir
d'un bord conforme $(X,[\gamma])$, Fefferman et Graham \cite{FefGra85}
ont trouv{\'e} un d{\'e}veloppement formel pour une m{\'e}trique d'Einstein
asymptotiquement hyperbolique $g$ sur un voisinage tubulaire
$M=]0,\epsilon[× X$. En coordonn{\'e}es g{\'e}od{\'e}siques, la m{\'e}trique s'{\'e}crit
$g=\frac{dx^2+g_x}{x^2}$, o{\`u} $g_x$ est une m{\'e}trique sur $\{x\}× X$,
admettant un d{\'e}veloppement
\begin{equation}
  \label{eq:2}g_x = 
  \begin{cases}
    \gamma + x^2 g_2 + x^3 g_3 + \cdots + x^n g_n + \cdots ,
      \quad n \text{ impair} \\
    \gamma + x^2 g_2 + x^3 g_3 + \cdots + x^n \ln x \, G + x^n g_n +
    \cdots,
      \quad n \text{ pair}
  \end{cases}
\end{equation}
o{\`u} tous les termes du d{\'e}veloppement sont formellement d{\'e}termin{\'e}s, {\`a}
l'exception du terme $g_n$. L'existence g{\'e}n{\'e}rale d'un d{\'e}veloppement
polyhomog{\`e}ne est prouv{\'e}e dans \cite{ChrDelLeeSki05,Hel}. Il en r{\'e}sulte que
deux m{\'e}triques ayant m{\^e}me infini conforme et m{\^e}me terme ind{\'e}termin{\'e},
{\'e}crites dans des coordonn{\'e}es g{\'e}od{\'e}siques, co{\"\i}ncident {\`a} un ordre infini
pr{\`e}s de $X$. Nous prouvons qu'elles sont {\'e}gales :

\begin{thm}\label{th:1}
  Soient deux m{\'e}triques d'Einstein asymptotiquement hyperboliques
  ayant le m{\^e}me infini conforme et le m{\^e}me terme ind{\'e}termin{\'e} ; alors,
  {\`a} un diff{\'e}omorphisme {\'e}gal {\`a} l'identit{\'e} au bord pr{\`e}s, les deux
  m{\'e}triques sont {\'e}gales dans un voisinage du bord.
\end{thm}

Bien entendu, il y a un analogue infinit{\'e}simal. Si on a une m{\'e}trique
d'Einstein asymptotiquement hyperbolique $g$, alors une d{\'e}formation
infinit{\'e}simale d'Einstein est un 2-tenseur sym{\'e}trique $\dg$,
satisfaisant la lin{\'e}arisation de l'{\'e}quation $\Ric(g)=-ng$, {\`a} savoir
\begin{equation}
  \label{eq:24}
  \frac12 \big(\nabla^*\nabla\dg + \Ric\circ\dg + \dg\circ\Ric -
  2\overset{\circ}{R}\dg \big) -\delta^*\delta\dg - \frac12 \delta^*d\Tr\dg = - n
  \dg .
\end{equation}
{\'E}crivons $g$ en jauge g{\'e}od{\'e}sique, comme ci-avant.  En agissant par un
diff{\'e}omorphisme infinit{\'e}simal, on peut mettre aussi $\dg$ dans une
jauge g{\'e}od{\'e}sique ($\partial_x\lrcorner \dg=0$) ; dans une telle jauge,
$\dg$ admet un d{\'e}veloppement similaire {\`a} (\ref{eq:2}). En particulier,
si la d{\'e}formation infinit{\'e}simale $\dg$ ne modifie pas la m{\'e}trique
conforme au bord, alors, toujours en jauge g{\'e}od{\'e}sique, on a
\begin{equation}
  \label{eq:25}
  \dg = \frac1{x^2} ( x^n \dg_n + \cdots ),
\end{equation}
o{\`u} tout le d{\'e}veloppement est d{\'e}termin{\'e} par le seul terme ind{\'e}termin{\'e}
$\dg_n$, un 2-tenseur sym{\'e}trique sur le bord (il est de plus {\`a}
divergence et trace nulles). On peut alors {\'e}noncer :

\begin{thm}\label{th:2}
  Soit $g$ une m{\'e}trique d'Einstein asymptotiquement hyperbolique, et
  $\dg$ une d{\'e}formation infinit{\'e}simale d'Einstein de $g$ d{\'e}finie pr{\`e}s
  du bord ; {\'e}crivons $\dg$ en jauge g{\'e}od{\'e}sique : si $|\dg|_g=o(x^n)$
  alors $\dg=0$ (toujours dans la jauge g{\'e}od{\'e}sique).
\end{thm}

Une version locale est vraie :

\begin{thm}\label{th:3}
  Les th{\'e}or{\`e}mes \ref{th:1} et \ref{th:2} restent vrais, mutatis
  mutandis, pour des m{\'e}triques d{\'e}finies seulement localement pr{\`e}s d'un
  point du bord.
\end{thm}

Un corollaire de la d{\'e}monstration, sugg{\'e}r{\'e} par M.~Herzlich, est
le cas plus facile, mais apparemment ouvert, d'un bord {\`a} distance finie :

\begin{thm}\label{th:4}
  Deux m{\'e}triques d'Einstein sur une vari{\'e}t{\'e} {\`a} bord $M$, de m{\^e}me
  constante d'Einstein, dont la restriction au bord $X=\partial M$ et la
  seconde forme fondamentale sur $X$ co{\"\i}ncident, diff{\`e}rent pr{\`e}s du
  bord par un diff{\'e}omorphisme induisant l'identit{\'e} sur $X$.
\end{thm}

Bien s{\^u}r, la version infinit{\'e}simale du th{\'e}or{\`e}me \ref{th:4}, que l'on
n'{\'e}noncera pas, reste valable aussi pour un bord {\`a} distance finie. En
revanche, une version locale de ce th{\'e}or{\`e}me semble plus difficile {\`a}
obtenir par les techniques d{\'e}velopp{\'e}es dans cet article.

La difficult{\'e} habituellement rencontr{\'e}e avec les {\'e}quations d'Einstein
est le manque d'ellipticit{\'e} de l'{\'e}quation sur la m{\'e}trique, {\'e}crite dans
des coordonn{\'e}es g{\'e}od{\'e}siques. On y rem{\'e}die en observant que, en
coordonn{\'e}es g{\'e}od{\'e}siques, la seconde forme fondamentale des tranches
satisfait une {\'e}quation elliptique non lin{\'e}aire, d{\'e}pendant de la
m{\'e}trique, avec second membre impliquant la courbure moyenne et ses
d{\'e}riv{\'e}es. Il y a alors une comp{\'e}tition entre les estimations de
Carleman adapt{\'e}es au cas asymptotiquement hyperbolique \cite{Maz91c},
et les termes d'erreur provenant de la m{\'e}trique et de la courbure
moyenne, mais le contr{\^o}le de ceux-ci par la seconde forme
fondamentale, gr{\^a}ce {\`a} la partie radiale de l'{\'e}quation d'Einstein,
s'av{\`e}re suffisant.

Dans cet article, on se limite au cas o{\`u} la m{\'e}trique au bord est de
classe $C^\infty$. Les estimations de Carleman n'utilisent que deux
d{\'e}riv{\'e}es de la m{\'e}trique, donc il est clair que nos r{\'e}sultats d'unique
continuation restent valables pour des m{\'e}triques au bord de r{\'e}gularit{\'e}
finie (par exemple $C^{3,\alpha}$), sous r{\'e}serve de supposer a priori
que les deux m{\'e}triques co{\"\i}ncident {\`a} un ordre infini pr{\`e}s du bord. En
r{\'e}alit{\'e}, la question naturelle consiste {\`a} supposer seulement une
co{\"\i}ncidence {\`a} l'ordre $o(x^n)$, et une r{\'e}gularit{\'e} $C^{2,\alpha}$,
voire m{\^e}me $C^\alpha$, sur le bord. La question d'unique continuation
pose alors des probl{\`e}mes plus techniques, qui seront abord{\'e}s ailleurs.
Il reste aussi {\`a} explorer la g{\'e}n{\'e}ralisation de ces r{\'e}sultats {\`a}
d'autres comportements asymptotiques, comme ceux {\'e}tudi{\'e}s dans
\cite{Biq00}.

\noindent\emph{Remerciements.} Luc Robbiano a r{\'e}pondu {\`a} mes questions
sur la continuation unique {\`a} un stade pr{\'e}liminaire de ce travail.
J'ai b{\'e}n{\'e}fici{\'e} de nombreuses discussions avec Marc Herzlich, qui m'a
indiqu{\'e} que M.~Anderson et lui travaillent sur une approche diff{\'e}rente
de la continuation unique pour les m{\'e}triques d'Einstein. Enfin, les
commentaires de Rafe Mazzeo furent tr{\`e}s utiles.

\section{Les {\'e}quations}
\label{sec:les-equations}

Nous utiliserons la coordonn{\'e}e $r=-\ln x$, de sorte que le mod{\`e}le
$\frac{dx^2+\gamma}{x^2}$ devient $dr^2+e^{2r}\gamma$. Soit une m{\'e}trique en
coordonn{\'e}es g{\'e}od{\'e}siques, $g=dr^2+\bg$. La seconde forme fondamentale
des tranches $\{x\}× X$ est $\II=-\frac12 \LL_{\partial_r}\bg$, dont on
notera la trace $H=\Tr^{\bg}(\II)$. Pour une m{\'e}trique asymptotiquement
hyperbolique, on a $\II \to -1$ et $H \to -n$. L'{\'e}quation $\Ric^g=-ng$
se d{\'e}compose alors en le syst{\`e}me, sur chaque tranche $\{x\}× X$,
\begin{align}
  \overline{\Ric} - H \II + \LL_{\partial_r}\II + 2 \II^2 &= -n \bg \label{eq:3}\\
  \bar{\delta}\II + dH &= 0 \label{eq:4}\\
  \partial_rH - |\II|^2 &= -n .\label{eq:5}
\end{align}
Ici $\overline{\Ric}$, $\bar{\delta}$,\ldots se rapportent {\`a} la m{\'e}trique
$\bg$, et $\II^2$, vu comme endomorphisme sym{\'e}trique pour $\bg$, est
le carr{\'e} de $\II$. Par la formule (\ref{eq:24}) pour la lin{\'e}arisation
du tenseur de Ricci, on a 
\begin{equation}
\begin{split}
 \LL_{\partial_r}\bRic&=\mathrm{d}_{\bg}\bRic(\LL_{\partial_r}\bg)
                =-2\mathrm{d}_{\bg}\bRic(\II)\\
 &=-2\big\{ \frac12\big( \bn^*\bn\II + \bRic\circ\II + \II\circ\bRic -
                        2\overset{\circ}{\bar R}\II \big)
            -\bar{\delta}^*\bar{\delta}\II-\frac12
            \bar{\delta}^*d\Tr^{\bg}(\II)\big\} .
\end{split}\label{eq:26}
\end{equation}
En appliquant $\LL_{\partial_r}$ {\`a} l'{\'e}quation (\ref{eq:3}), et en
utilisant (\ref{eq:4}) qui s'apparente {\`a} une condition de jauge sur
$\II$, on d{\'e}duit de (\ref{eq:26}) l'{\'e}quation
\begin{multline}
  \label{eq:7}
  -\bn^*\bn\II+\LL_{\partial_r}^2\II-H\LL_{\partial_r}\II - (\partial_rH)\II+2(\II\circ\LL_{\partial_r}\II+2\II\circ\II\circ\II+\LL_{\partial_r}\II\circ\II) \\
- \bRic\circ\II - \II\circ\bRic + 2\overset{\circ}{\bar R}\II-2n\II -
\bar{\delta}^*dH = 0 .
\end{multline}
On remarquera qu'en faisant dispara{\^\i}tre $H$ de cette {\'e}quation gr{\^a}ce {\`a}
l'{\'e}quation (\ref{eq:4}), on trouverait une autre {\'e}quation, non
elliptique, sur $\II$. C'est pourquoi on pr{\'e}f{\`e}re ici traiter les
termes en $H$ comme un second membre.

\section{Estimation de $g$ et $H$}
\label{sec:estimation-de-g}

Pour une fonction $f$ sur $[r_0,+\infty[$, telle que $|f|=O(e^{-\lambda r})$
pour tout $\lambda>0$, une int{\'e}gration par parties donne imm{\'e}diatement
\begin{equation}
  \label{eq:8}
  \int_{r_0}^{+\infty} |\partial_rf|^2 e^{2\lambda r} dr \geq \lambda^2 \int_{r_0}^{+\infty} |f|^2
  e^{2\lambda r} dr .
\end{equation}
Aucune condition sur $f$ en $r_0$ n'est n{\'e}cessaire.

Supposons maintenant que deux m{\'e}triques d'Einstein asymptotiquement
hyperboliques $g$ et $g_0$ aient m{\^e}mes infini conforme et terme
ind{\'e}termin{\'e}. Une m{\'e}trique, choisie dans la classe conforme, d{\'e}termine
des coordonn{\'e}es g{\'e}od{\'e}siques dans un voisinage du bord \cite{Gra00},
dans lesquelles $g$ et $g_0$ co{\"\i}ncident donc {\`a} un ordre infini au
bord.  Nous pouvons nous restreindre {\`a} un voisinage $]r_0,+\infty[×X$
du bord sur lequel $g$ et $g_0$ sont mutuellement born{\'e}es, et toutes
les d{\'e}riv{\'e}es de $g-g_0$ restent born{\'e}es par rapport {\`a} $g_0$. De
l'estimation (\ref{eq:8}), appliqu{\'e}e dans des syst{\`e}mes de coordonn{\'e}es,
r{\'e}sulte alors imm{\'e}diatement, pour $\lambda>2$,
\begin{equation}
  \label{eq:6}
  \int_{r_0}^{+\infty} |\II-\II_0|_{g_0}^2 e^{2\lambda r} dr \geq C^{-1} \lambda^2
  \int_{r_0}^{+\infty} |g-g_0|_{g_0}^2 e^{2\lambda r} dr
\end{equation}
pour une constante $C$ ind{\'e}pendante de $\lambda$. Dans la suite, la
constante $C$ pourra grandir {\`a} chaque ligne, mais restera ind{\'e}pendante
de $\lambda$. On peut aussi d{\'e}river l'{\'e}quation
$\II-\II_0=-\frac12 \LL_{\partial_r}(g-g_0)$ : comme toutes les courbures et
secondes formes fondamentales sont born{\'e}es, les commutations
n'introduisent pas de nouveau terme et on obtient
\begin{equation}
  \label{eq:10}
  \int_{r_0}^{+\infty} \sum_0^k |\nabla_0^j(\II-\II_0)|_{g_0}^2 e^{2\lambda r} dr \geq C^{-1} \lambda^2
  \int_{r_0}^{+\infty} \sum_0^k |\nabla_0^j(g-g_0)|_{g_0}^2 e^{2\lambda r} dr.
\end{equation}

De m{\^e}me, de l'{\'e}quation (\ref{eq:5}) on d{\'e}duit
$\partial_r(H-H_0)=|\II|_g^2-|\II_0|_{g_0}^2$ ; puisque $g$ et toutes ses
d{\'e}riv{\'e}es demeurent uniform{\'e}ment born{\'e}es par rapport {\`a} $g_0$ dans le
voisinage consid{\'e}r{\'e}, on a $\big| |\II|_g^2-|\II_0|_{g_0}^2 \big| \leq C
(|\II-\II_0|_{g_0}+|g-g_0|_{g_0})$, et de (\ref{eq:8}) et (\ref{eq:6})
r{\'e}sulte l'estimation
\begin{equation}
  \label{eq:9}
  \int_{r_0}^{+\infty}|\II-\II_0|_{g_0}^2 e^{2\lambda r} dr \geq C^{-1} \lambda^2
  \int_{r_0}^{+\infty}|H-H_0|^2 e^{2\lambda r} dr .
\end{equation}
En prenant des d{\'e}riv{\'e}es de (\ref{eq:5}), comme ci-dessus,
la m{\^e}me estimation reste valable pour les d{\'e}riv{\'e}es :
\begin{equation}
  \label{eq:11}
  \int_{r_0}^{+\infty}\sum_0^k |\nabla_0^j(\II-\II_0)|_{g_0}^2 e^{2\lambda r} dr \geq C^{-1} \lambda^2
  \int_{r_0}^{+\infty}\sum_0^k |\nabla_0^j(H-H_0)|^2 e^{2\lambda r} dr .
\end{equation}

\section{Estimation de Carleman et preuve des th{\'e}or{\`e}mes \ref{th:1} et \ref{th:2}}
\label{sec:estim-de-carl}

Rappelons ici la version de Mazzeo de l'estimation de Carleman
\cite[th{\'e}or{\`e}me 7]{Maz91c} : si $P$ est un op{\'e}rateur elliptique pour
une m{\'e}trique asymptotiquement hyperbolique $g$, de la forme
\begin{equation}
Ps=\nabla^*\nabla s+\text{ termes d'ordre inf{\'e}rieur born{\'e}s en }\nabla s \text{ et }s,\label{eq:13}
\end{equation}
alors il existe $r_0>0$ tel que pour tous $\lambda\gg 0$ et $s$ {\`a} support
dans $\{r> r_0\}$ on ait
\begin{equation}
  \label{eq:12}
  \int_{r\geq r_0} |Ps|^2 e^{2\lambda r} \vol^g \geq
  C^{-1} \int_{r\geq r_0} \big( \lambda^{-1}|\nabla^2s|^2+\lambda|\nabla s|^2+\lambda^3|s|^2 \big)
    e^{2\lambda r} \vol^g .
\end{equation}
En fait le terme $|\nabla^2s|^2$ n'est pas {\'e}crit dans l'{\'e}nonc{\'e} cit{\'e}, mais
il est {\'e}vident de revoir la d{\'e}monstration pour montrer sa pr{\'e}sence avec
le coefficient indiqu{\'e}. En outre, les constantes $r_0$ et $C$ peuvent
{\^e}tre choisies uniformes si $g$ et les coefficients de $P$ varient de
mani{\`e}re born{\'e}e.

Appliquons {\`a} pr{\'e}sent cette estimation {\`a} la d{\'e}monstration du th{\'e}or{\`e}me
\ref{th:1}. Comme pr{\'e}c{\'e}demment, nous avons donc deux m{\'e}triques
d'Einstein $g$ et $g_0$ qui co{\"\i}ncident {\`a} un ordre infini dans des
coordonn{\'e}es g{\'e}od{\'e}siques, et nous analysons la diff{\'e}rence des {\'e}quations
(\ref{eq:7}) pour les m{\'e}triques $g$ et $g_0$. Par exemple, le premier
terme $\bn^*\bn\II-\bn_0^*\bn_0\II_0$ se d{\'e}compose en
\begin{equation}
  \label{eq:14}
  \bn_0^*\bn_0(\II-\II_0) + (\bn^*\bn-\bn_0^*\bn_0)\II ;
\end{equation}
dans cette expression, le premier terme est la lin{\'e}arisation en $g_0$
et, puisque toutes les d{\'e}riv{\'e}es de $g$, et donc de $\II$ aussi, sont
suppos{\'e}es born{\'e}es, le second terme est contr{\^o}l{\'e} en chaque point par
$C(|g-g_0|+|\bn_0(g-g_0)|)$.

Plus g{\'e}n{\'e}ralement, notons $\Phi(g,H)(\II)$ le membre de gauche de
l'{\'e}quation (\ref{eq:7}), et $P=\frac{\partial\Phi}{\partial\II}(g_0,H_0)$ la
lin{\'e}arisation de l'{\'e}quation en $\II_0$. Alors $\Phi(g_0,H_0)(\II_0)=0$ et
\begin{equation}
  \label{eq:15}
  \Phi(g,H)(\II) =  P(\II-\II_0) + Q(g,H)(\II),
\end{equation}
et une lecture attentive de chaque terme comme dans l'{\'e}quation
(\ref{eq:14}) donne l'estimation
\begin{multline}
  \label{eq:16}
  |Q(g,H)(\II)| \leq \epsilon \big( |\II-\II_0| + |\nabla_0(\II-\II_0)| \big) \\
    + C \big( |g-g_0|+|\nabla_0(g-g_0)|+|\nabla_0^2(g-g_0)| \\
    + |H-H_0| + |\nabla_0(H-H_0)| + |\nabla_0^2(H-H_0)| \big).
\end{multline}
La constante $\epsilon$ peut {\^e}tre prise petite, mais ce fait ne sera
pas utilis{\'e} dans la suite. Combinant avec les estimations
(\ref{eq:10}) et (\ref{eq:11}), on d{\'e}duit l'estimation
\begin{multline}
  \label{eq:17}
  \int_{r\geq r_0}|Q(g,h)(\II)|^2 e^{2\lambda r}\vol \\ \leq
  \int_{r\geq r_0}\big( \epsilon( |\II-\II_0|^2 + |\nabla_0(\II-\II_0)|^2 )
  + C \lambda^{-2} \sum_0^2 |\nabla_0^j(\II-\II_0)|^2 \big) e^{2\lambda r} \vol
\end{multline}
D'un autre c{\^o}t{\'e}, l'op{\'e}rateur $P$ est un op{\'e}rateur du second ordre,
dont les termes de plus haut degr{\'e} sont $-\bn_0^*\bn_0+\LL_{\partial_r}^2$,
et les autres termes sont uniform{\'e}ment born{\'e}s, donc il satisfait
l'estimation de Carleman (\ref{eq:12}), que nous pouvons appliquer
apr{\`e}s avoir coup{\'e} $\II-\II_0$ entre $r_0$ et $r_0+1$, d'o{\`u} :
\begin{multline}
  \label{eq:18}
  \int_{r\geq r_0} |P(\II-\II_0)|^2 e^{2\lambda r} \vol
 + \int_{r_0\leq r\leq r_0+1} \sum_0^2 |\nabla_0^j(\II-\II_0)|^2 e^{2\lambda r} \vol\\
 \geq C^{-1} \int_{r\geq r_0}
 (\lambda^{-1}|\nabla_0^2(\II-\II_0)|^2+\lambda|\nabla_0(\II-\II_0)|^2+\lambda^3|\II-\II_0|^2)
 e^{2\lambda r}\vol . 
\end{multline}

L'{\'e}quation $P(\II-\II_0)+Q(g,H)(\II)=0$, avec les deux estimations
pr{\'e}c{\'e}dentes, fournit finalement pour $\lambda\gg 0$,
\begin{multline}
  \label{eq:19}
  \int_{r_0\leq r\leq r_0+1} \sum_0^2 |\nabla_0^j(\II-\II_0)|^2 e^{2\lambda r} \vol\\
 \geq C^{-1} \int_{r\geq r_0+1}
 (\lambda^{-1}|\nabla_0^2(\II-\II_0)|^2+\lambda|\nabla_0(\II-\II_0)|^2+\lambda^3|\II-\II_0|^2)
 e^{2\lambda r}\vol.
\end{multline}
En faisant $\lambda\to+\infty$, on voit qu'il faut $\II-\II_0=0$ sur $r\geq
r_0+1$, donc $g=g_0$.

Le th{\'e}or{\`e}me \ref{th:2} est la version infinit{\'e}simale du th{\'e}or{\`e}me
\ref{th:1}, sa d{\'e}monstration est plus facile : on {\'e}crit la
lin{\'e}arisation du syst{\`e}me (\ref{eq:3}), (\ref{eq:4}) et (\ref{eq:5}),
que l'on traite de mani{\`e}re similaire.

\section{Variante locale pour un bord {\`a} l'infini}
\label{sec:variante-locale-pour}

Pour d{\'e}montrer le th{\'e}or{\`e}me \ref{th:3} dans le cas asymptotiquement
hyperbolique, on localise le probl{\`e}me au voisinage d'un point du bord
en se rappelant du mod{\`e}le $g_{n+1}=dr^2+\cosh^2(r)g_n$ pour l'espace
hyperbolique r{\'e}el. L'infini conforme quand $r\to+\infty$ est alors juste
une boule hyperbolique.

Le r{\'e}sultat local de continuation unique pour un op{\'e}rateur elliptique
dans \cite{Maz91c} est montr{\'e} en utilisant des coordonn{\'e}es obtenues
par les g{\'e}od{\'e}siques issues d'une hypersurface
$x^2+x_1^2+\cdots+x_n^2=\epsilon$, ce qui donne bien un mod{\`e}le du type
indiqu{\'e} avant. Cependant, nous ne pouvons pas proc{\'e}der ainsi, puisque
l'hypoth{\`e}se de co{\"\i}ncidence des m{\'e}triques sur le bord $x=0$ serait
alors perdue.  Il faut donc travailler {\`a} partir du bord.

Soit une m{\'e}trique d'Einstein asymptotiquement hyperbolique $g$,
d'infini conforme $[\gamma]$. Pr{\`e}s d'un point $p\in X=\partial M$, on fixe une
m{\'e}trique $\gamma$ dans la classe conforme, et une petite boule $B\subset X$ pour
$\gamma$, centr{\'e}e en $p$. Soit $y$ la distance au bord $\partial B$ pour $\gamma$,
prolong{\'e}e de mani{\`e}re lisse et strictement positive {\`a} l'int{\'e}rieur de
$B$. Pr{\`e}s de $X\subset M$ on a donc $g=\frac{dx^2+\gamma}{x^2}$, et dans $X=\{x=0\}$,
pr{\`e}s de $\partial B\subset X$, on a de plus la d{\'e}composition $\gamma|_X=dy^2+\bga$, o{\`u}
$\bga$ est une famille de m{\'e}triques sur $\partial B$, param{\'e}tr{\'e}e par $y$. En
prenant alors sur $M$ des coordonn{\'e}es $(x,y,x_2,\dots,x_n)$, o{\`u}
$(x_2,\dots,x_n)$ sont des coordonn{\'e}es sur $\partial B$, on obtient donc
\begin{equation}
  \label{eq:35}
  \gamma=dy^2+\bga+\gamma', \quad \bga=\bga_{ij}(y,x_2,\dots,x_n)dx_idx_j,\quad \gamma'=O(x) .
\end{equation}
On effectue alors le changement de variables
\begin{equation}
  \label{eq:29}
  x=v \frac{2u_0}{1+u_0^2}, \quad y=v \frac{1-u_0^2}{1+u_0^2},
\end{equation}
de sorte que
\begin{equation}
  \label{eq:28}
  g = \frac{dx^2+dy^2+\bga+\gamma'}{x^2} =
  \frac{1}{u_0^2}\big(du_0^2+\frac{(1+u_0^2)^2}4 \frac{dv^2+\bga+\gamma'}{v^2}\big) .
\end{equation}
Ce mod{\`e}le nous incite {\`a} choisir dans la classe conforme de $\gamma$ sur
$B$ la m{\'e}trique asymptotiquement hyperbolique
\begin{equation}
  \label{eq:30}
  \tilde \gamma = \frac{\gamma}{4y^2} ,
\end{equation}
puis {\`a} {\'e}crire la m{\'e}trique dans les coordonn{\'e}es g{\'e}od{\'e}siques
correspondantes : 
\begin{equation}
g=\frac{du^2+\bg}{u^2},
\label{eq:20}
\end{equation}
avec $u=u_0e^\phi$ une « {\'e}quation sp{\'e}ciale » de $X$, c'est-{\`a}-dire
satisfaisant la condition
\begin{equation}
  \label{eq:31}
  |du|_{u^2 g} = 1, \quad u^2 g|_X = \tilde \gamma .
\end{equation}
Graham \cite{Gra00} montre que cette {\'e}quation aux d{\'e}riv{\'e}es partielles
du premier ordre, non caract{\'e}ristique le long de $X$, a toujours une
solution dans un voisinage de $X$, mais ici la m{\'e}trique $\tilde \gamma$
est singuli{\`e}re sur $\partial B$ et il faut {\^e}tre plus prudent. L'{\'e}quation se
r{\'e}crit
\begin{equation}
  \label{eq:33}
  (\nabla^{u_0^2g}u_0)\cdot \phi + \frac {u_0}2 |d\phi|^2_{u_0^2g} =
  \frac{1-|du_0|^2_{u_0^2g}}{2u_0}, \quad \phi|_X = 0 .
\end{equation}
{\`A} partir de (\ref{eq:28}), on voit que le second membre est
$O(v^2u_0)$, si bien que l'existence d'une solution dans un voisinage
de $\{u_0=0\}$ ne pose pas de difficult{\'e}, et celle-ci s'annule le long
de $\{v=0\}$.  Dans ces nouvelles coordonn{\'e}es $(u,v,x_2,\dots,x_n)$,
la m{\'e}trique $g$ a donc la forme (\ref{eq:20}), et chaque m{\'e}trique
$\bg$ est elle-m{\^e}me asymptotiquement hyperbolique sur $B$, avec infini
conforme ind{\'e}pendant de $u$ et {\'e}gal {\`a} $\bga$.

Si deux telles m{\'e}triques $g_1$ et $g_2$ co{\"\i}ncident {\`a} un ordre infini
le long de $X$, il est clair aussi que les deux fonctions
correspondantes $\phi_1$ et $\phi_2$, et les m{\'e}triques $\bg_1$ et $\bg_2$,
co{\"\i}ncident {\`a} un ordre infini le long de $\{x=0\}=\{u=0\} \cup \{v=0\}$.
On peut alors faire les int{\'e}grations par parties menant {\`a} l'estimation
de Carleman (\ref{eq:12}), et le reste de la d{\'e}monstration est
inchang{\'e}e. Le th{\'e}or{\`e}me \ref{th:3} est donc d{\'e}montr{\'e}.

\section{D{\'e}monstration du th{\'e}or{\`e}me \ref{th:4}}
\label{sec:demonstr-des-theor}

Le cas d'un bord {\`a} distance finie se traite de la fa{\c c}on suivante.
Soient $g$ et $g_0$ deux m{\'e}triques lisses sur la vari{\'e}t{\'e} {\`a} bord $M$.
On identifie {\`a} nouveau un voisinage de $X$ dans $M$ au produit
$[0,\epsilon[×X$, avec coordonn{\'e}es g{\'e}od{\'e}siques, de sorte que $g=dx^2+\bg$ et
$g_0=dx^2+\bg_0$. Le syst{\`e}me des {\'e}quations d'Einstein (\ref{eq:3}),
(\ref{eq:4}) et (\ref{eq:5}), apr{\`e}s avoir modifi{\'e} la constante
d'Einstein et remplac{\'e} $r$ par $x$, dit imm{\'e}diatement que toutes les
d{\'e}riv{\'e}es radiales de $\bg$ et $\bg_0$ co{\"\i}ncident en $x=0$. On applique
alors la m{\^e}me m{\'e}thode que pour le th{\'e}or{\`e}me \ref{th:1}, en faisant les
estimations par rapport {\`a} la m{\'e}trique lisse $g_0$ : ainsi l'estimation
(\ref{eq:8}) devient-elle, d{\`e}s que $\lambda>-1/2$, pour une fonction $f$
s'annulant {\`a} un ordre infini en $x=0$,
\begin{equation}
  \label{eq:21}
  \int_0^{x_0} |\partial_x f|^2 x^{-2\lambda} dx
 \geq \lambda^2 \int_0^{x_0} |f|^2 x^{-2\lambda-2} dx .
\end{equation}
Plus g{\'e}n{\'e}ralement, (\ref{eq:10}) devient
\begin{multline}
  \label{eq:22}
  \int_0^{x_0} |\nabla_0^2(\II-\II_0)|^2 x^{-2\lambda} dx \geq \\
 C^{-1} \lambda^2 \int_0^{x_0} \big( |\nabla_0^2(g-g_0)|^2 + \lambda^2 x^{-2}
 |\nabla_0(g-g_0)|^2 + \lambda^4 x^{-4} |g-g_0|^2 \big) x^{-2\lambda} dx
\end{multline}
et on a le m{\^e}me contr{\^o}le sur $H-H_0$ {\`a} la place de (\ref{eq:11}).
Enfin, l'estimation de Carleman (\ref{eq:12}) devient, pour $\lambda\gg 0$
et $x_0$ assez petit,
\begin{multline}
  \label{eq:23}
  \int_0^{x_0} |Ps|^2 x^{-2\lambda} \vol^{g_0} \geq \\
 C^{-1} \lambda^{-1} \int_0^{x_0} \big( |\nabla_0^2s|^2 + \lambda^2x^{-2}|\nabla_0s|^2
 + \lambda^4x^{-4}|s|^2 \big) x^{-2\lambda} \vol^{g_0} .
\end{multline}
{\`A} partir de ces estimations, \emph{meilleures} que dans le cas
asymptotiquement hyperbolique, il devient clair que la d{\'e}monstration
du th{\'e}or{\`e}me \ref{th:1} s'applique encore.

Bien entendu, le m{\^e}me raisonnement s'applique au cas infinit{\'e}simal. En
revanche, la m{\'e}thode de localisation employ{\'e}e dans la section
\ref{sec:variante-locale-pour} dans le cas asymptotiquement
hyperbolique ne peut plus {\^e}tre appliqu{\'e}e ici, puisque c'est toute la
m{\'e}trique, plut{\^o}t que sa classe conforme, qui est fix{\'e}e au bord.

\begin{rem}[Unique continuation en un point]
  La m{\^e}me d{\'e}monstration est valable pour deux m{\'e}triques d{\'e}finies sur
  un ouvert au voisinage d'un point $p$, avec $x$ d{\'e}signant la
  distance {\`a} $p$ : il suffit de raisonner en coordonn{\'e}es normales. On
  obtient ainsi, sans recours aux coordonn{\'e}es harmoniques, une
  d{\'e}monstration alternative du fait que deux m{\'e}triques d'Einstein, de
  m{\^e}me constante d'Einstein, {\'e}gales au point $p$, dont les courbures
  et toutes leurs d{\'e}riv{\'e}es covariantes co{\"\i}ncident au point $p$, sont
  {\'e}gales.
\end{rem}

\bibliographystyle{smfplain}
\bibliography{biblio,biquard,uc}

\end{document}